\documentclass{article}
\usepackage{amsmath,theorem,latexsym,amssymb}

\newtheorem{theorem}{Theorem}[section]

\begin{document}

\title{Non-divergence harmonic maps}


%\address{Max Planck Institute for Mathematics in the Sciences, Inselstrasse 22, 04103 Leipzig,
%Germany}
%\curraddr{Department of Mathematics and Statistics,
%Case Western Reserve University, Cleveland, Ohio 43403}
%\email{jost@mis.mpg.de}

\author{J\"urgen Jost and Fatma Muazzez \c Sim\c sir}
%\address{\c{C}\i nar Sokak No: 115/2 06170 Ankara, Turkey}
%\email{fsimsir@gmail.com}

%    General info
%\subjclass{Primary 53B05; Secondary 58E20}
%\date{January 1, 1994 and, in revised form, June 22, 1994.}

%\keywords{Affine flat connection, affine harmonic map, non-divergence type
%invariant elliptic operator}

\maketitle
\begin{abstract}
We describe work on solutions of certain non-divergence type and
therefore non-variational  elliptic and
parabolic systems on manifolds. These systems include Hermitian and
affine harmonics which should become useful tools for studying Hermitian and
 affine manifolds, resp. A key point is that in addition to
the standard condition of nonpositive image curvature that is well
known and understood in the theory of ordinary harmonic maps (which
arise from a variational problem), here we also need in addition a
global topological condition to guarantee the existence of
solutions.
\end{abstract}

\section*{Introduction}
In this paper, we shall describe concepts and tools from geometric
analysis that we have developed for particular classes of manifolds,
namely Hermitian and affine ones. Hermitian manifolds are complex
manifolds that are also equipped with an Hermitian metric. Similarly,
an affine manifold can be equipped with a Riemannian metric as an
auxiliary structure.  
Here, a manifold is said to be flat or affine if it admits an atlas whose coordinate changes
are affine transformations.

Basic tools of Riemannian geometry are the geodesics and
their higher dimensional generalizations, the harmonic maps. They are
the critical points of an energy integral that involves the metric. Therefore, they are backed
by a variational structure. This depends on the Levi-Civit\`a
connection underlying the Riemannian metric. A Hermitian manifold,
however, naturally possesses a different connection, the complex one
that respects the complex structure. This connection is different from
the Levi-Civit\`a connection unless the manifold is
K\"ahler. Similarly, an affine manifold carries a flat affine
connection that has nothing to do with the Levi-Civit\`a connection of
the auxiliary Riemannian metric. In particular, that Riemannian metric
need not be flat.\\
Thus, harmonic maps are not naturally defined on such manifolds, and
the main point of this paper is to discuss suitable substitutes. Thus, Hermitian harmonic maps, as
introduced and studied in \cite{JY}, are
defined through the complex connection, and affine harmonic
maps, as introduced and studied in \cite{JS}, are determined by the  affine connection, 
and the resulting equations do not 
satisfy a variational principle. This is already the case for affine
geodesics, as there is in general no length or energy functional that
they could locally minimize. Also, the Euler-Lagrange equations of
variational problems necessarily have a special, divergence-type
structure which in general the affine harmonic map equations do not
possess. The absence of a variational structure makes the analysis more difficult. 
Therefore, we need an additional global non-triviality condition to guarantee the
existence of an affine harmonic map in a given homotopy class. As in the case
of ordinary harmonic maps, nonpositive  curvature of the target manifold is also required. 

In this paper, we overview the results of \cite{JS}, and its connections with the previous work of 
\cite{JY}. 
For all geometric concepts and notations not explained here, as well
as for a recent treatment and survey of the theory of harmonic maps, we refer
to \cite{J3} as our standard reference. In particular, we shall use
the heat equation method as introduced in the seminal paper \cite{MR}
and applied in many subsequent papers in geometric analysis (see
\cite{J4} for a more detailed history). 
However, as we do not have a variational structure at our disposal,
 we cannot utilize the arguments of those papers and have to proceed
 rather differently.

%%%%%%%%%%%%%%%%%%%%%%%%%%%%%%%%%%%%%%%%%%%%%%%%%%%%%%%%%%%%%%%%%%%%%%%%%%%%%%%%%%%

\section{Coordinate transformations and invariant
  differential operators}
A Riemannian metric $\gamma$ on a manifold $M$ is locally, that is,
w.r.t. local coordinates $x^\alpha$,  of the form 
       \begin{equation}
        \label{100}
           \gamma = \gamma_{\alpha \beta}
                     dx^{\alpha} \otimes dx^{\beta}, 
       \end{equation} 
and under coordinate transformations $x=x(y)$, it transforms as
\begin{equation}
        \label{101}
          \gamma_{\alpha \beta}
                     dx^{\alpha} \otimes dx^{\beta}= \gamma_{\alpha
                       \beta}\frac{\partial x^\alpha}{\partial
                       y^\delta}\frac{\partial x^\beta}{\partial
                       y^\eta}dy^\delta dy^\eta =: h_{\delta \eta}dy^\delta dy^\eta.
       \end{equation} 
Therefore, the coefficients of the inverse metric tensor transform
according to
\begin{equation}
        \label{102}
          \gamma^{\alpha
                       \beta}\frac{\partial
                       y^\delta}{\partial x^\alpha}\frac{\partial
                       y^\eta}{\partial x^\beta}= h^{\delta \eta}.
       \end{equation} 
Now, the second derivative of a function $\phi$,
\begin{equation}
 \label{103}
\frac{\partial^2 \phi}{\partial x^\alpha \partial x^\beta}
\end{equation}
transforms into
\begin{equation}
 \label{104}
\frac{\partial^2 \phi}{\partial x^\alpha \partial x^\beta}\frac{\partial x^\alpha}{\partial
                       y^\delta}\frac{\partial x^\beta}{\partial
                       y^\eta}+ \frac{\partial \phi}{\partial
                       x^\alpha}\frac{\partial^2 x^\alpha}{\partial
                       y^\delta \partial y^\eta},
\end{equation}
that is, there is an additional term with {\it second} derivatives of the
coordinate transformation. Therefore, in general,
\begin{equation}
 \label{105}
\gamma^{\alpha \beta}\frac{\partial^2 \phi}{\partial x^\alpha \partial x^\beta}
\end{equation}
is not invariant -- unless that second derivative $\frac{\partial^2
  x^\alpha}{\partial y^\delta \partial y^\eta}$ vanishes for all
indices. This is the case when the coordinate transformations are
{\it affine linear}. In other words, on an affine manifold, the metric trace
of the second derivatives of a function $\phi$ is coordinate
invariant.

Analogously, on a complex manifold, we may consider a Hermitian metric
   \begin{equation}
        \label{106}
           \gamma = \gamma_{\alpha \bar{\beta}}
                     dz^{\alpha} \otimes dz^{\bar{\beta}}, 
       \end{equation} 
and the Hermitian trace 
\begin{equation}
 \label{107}
\gamma^{\alpha \bar{\beta}}\frac{\partial^2 \phi}{\partial z^\alpha \partial z^{\bar{\beta}}}
\end{equation}
is invariant under {\it holomorphic} coordinate transformations. 
More generally, the same applies for a map $f:M\to N$ from $M$ into
some Riemannian manifold $N$ in place of a function $\phi$. Thus, 
denoting the Christoffel symbols of $N$ in local coordinates by
$\Gamma^i_{j k}$, we have the invariant operator on an affine manifold $M$
\begin{equation}
      \label{108}
      \gamma^{\alpha {\beta}} ( \frac{\partial^2 f^i}{\partial x^\alpha \partial x^{{\beta}}}                  			 +\Gamma^{i}_{j k}\frac{\partial f^j}{\partial x^\alpha}
      \frac{\partial f^k}{\partial x^{{\beta}}} ), ~~i = 1, \ldots, n.
 \end{equation}
Similarly, on a complex manifold, we obtain the operator
  \begin{equation}
      \label{109}
      \gamma^{\alpha \bar{\beta}} ( \frac{\partial^2 f^i}{\partial z^\alpha \partial z^{\bar{\beta}}}                  			 +\Gamma^{i}_{j k}\frac{\partial f^j}{\partial z^\alpha}
      \frac{\partial f^k}{\partial z^{\bar{\beta}}} ), ~~i = 1,
      \ldots, n. 
 \end{equation}
We then call a solution $\phi$, resp., $f$ of 
\begin{equation}
 \label{110}
\gamma^{\alpha \beta}\frac{\partial^2 \phi}{\partial x^\alpha \partial
  x^\beta}=0; \quad \quad  ~ \gamma^{\alpha {\beta}} ( \frac{\partial^2 f^i}{\partial x^\alpha \partial x^{{\beta}}}                  			 +\Gamma^{i}_{j k}\frac{\partial f^j}{\partial x^\alpha}
      \frac{\partial f^k}{\partial x^{{\beta}}} ) = 0, ~~i = 1, \ldots, n
\end{equation} 
an affine harmonic function, resp., map. Analogously, a solution of
\begin{equation}
\label{111}
\gamma^{\alpha \bar{\beta}}\frac{\partial^2 \phi}{\partial z^\alpha
  \partial z^{\bar{\beta}}}=0; \quad \quad ~ \gamma^{\alpha \bar{\beta}} ( \frac{\partial^2 f^i}{\partial z^\alpha \partial z^{\bar{\beta}}}                  			 +\Gamma^{i}_{j k}\frac{\partial f^j}{\partial z^\alpha}
      \frac{\partial f^k}{\partial z^{\bar{\beta}}} )=0, ~~i = 1,
      \ldots, n
\end{equation}
on a complex manifold is called Hermitian harmonic. \\
We note that the equations (systems) (\ref{110}), (\ref{111}) are not in
divergence form, in contrast to the equation (system) for ordinary
harmonic functions (maps) on a Riemannian manifold. This makes the
existence and regularity theory more difficult. \\
In fact, ordinary harmonic functions (maps) on a Riemannian manifold
$M$ satisfy
\begin{eqnarray}
\label{112}
  \frac{1}{\sqrt{\det \gamma}}\frac{\partial}{\partial x^\alpha}(\sqrt{\det
  \gamma} \gamma^{\alpha \beta} \frac{\partial \phi}{\partial
  x^\beta})  &=& 0 \\
  \frac{1}{\sqrt{\det \gamma}}\frac{\partial}{\partial x^\alpha}(\sqrt{\det
\gamma} \gamma^{\alpha \beta}\frac{\partial f^i}{\partial
  x^\beta})+\gamma^{\alpha \beta}\Gamma^{i}_{j k}\frac{\partial f^j}{\partial x^\alpha}
      \frac{\partial f^k}{\partial x^{{\beta}}} &=& 0, ~~i = 1, \ldots, n,
\end{eqnarray} 
that is, some derivatives of the metric need to compensate second
derivatives of the function under coordinate changes in order to make the diferential equation
invariant. 

%%%%%%%%%%%%%%%%%%%%%%%%%%%%%%%%%%%%%%%%%%%%%%%%%%%%%%%%%%%%%%%%%%%%%%%%%%%%%%%%%%%%%%%%%%%%%%%

\section{Harmonic maps without variational or divergence structure}

A geometric structure usually induces a particular type of connection that preserves that structure.
In Riemannian geometry, the Levi-Civit\`a connection is the unique torsion free connection that preserves the
Riemannian metric. For a complex structure we get a canonical complex connection. Similarly, for 
an affine structure, we have the affine flat connection. As a result, different structures on the
same manifold induce different connections. To investigate such
structures, we need 
appropriate tools from geometric analysis. In Riemannian geometry,
geodesics and their higher dimensional analogues,
harmonic maps,  are such tools. Commonly, geodesics and harmonic maps
are defined in terms of a variational principle, as critical points of
the energy integral. This, however, is special for the Levi-Civit\'a
connection in Riemannian geometry and does not generalize to Hermitian
or affine geometry. Thus, as described above, we 
rather 
define such objects directly in terms of the relevant connection. We then
obtain an elliptic system that can be written in local coordinates,
but whose solutions are invariant under coordinate changes in the
appropriate category (differentiable, complex, affine). These
solutions then yield suitable classes of functions (when the
target is $\mathbb{R}$) or maps (when the target is a Riemannian
manifold). For instance, maps defined in this way, between Hermitian and Riemannian manifolds are called Hermitian harmonic  and  harmonic maps from affine flat to Riemannian manifolds are called affine harmonic.
The notion of Hermitian harmonic maps was first introduced and investigated by Jost and Yau \cite{JY} and that 
of affine harmonic maps was first introduced and investigated by Jost
and \c Sim\c sir \cite{JS}. In either case, a solution of the elliptic
system was obtained, under suitable conditions, from the associated
parabolic system.

Parabolic and elliptic systems with a nonlinearity as in the harmonic map problem  and
without a variational or divergence structure have been investigated by von Wahl \cite{W1}. 
However, he was mainly interested in boundary value problems on
Euclidean domains and  not
in the case of closed manifolds. Therefore, in order to treat the
central problem of analyzing  when the solution of the parabolic
system converges to that of elliptic one, 
a more global approach had to be  developed by Jost and Yau for the Hermitian harmonic maps in \cite{JY}. 
Extensions of existence and uniques results for the Dirichlet problem in the work of Jost and Yau
to noncompact but complete domain manifolds were first considered by Lei Ni \cite{N}. 
Subsequently, Grunau and K\"uhnel \cite{GK} developed a more flexible method. 
Throughout this work, harmonic map systems without a variational structure in which the underlying 
equations is of  non-divergence form will be called  non-divergence  harmonic maps.
                                   
%%%%%%%%%%%%%%%%%%%%%%%%%%%%%%%%%%%%%%%%%%%%%%%%%%%%%%%%%%%%%%%%%%%%%%%%%%%%%%%%%%%%%%%%%%%%%%%%%%%%%%%%%%%%%%%%%
		   
\subsection{Hermitian harmonic maps}

Let $M$ be a compact complex manifold with a Hermitian metric 
$(\gamma_{\alpha \bar{\beta}})_{\alpha, \beta = 1,\ldots, m}$ in local 
coordinates $z = (z^1, \ldots, z^m)$, and $N$ be a compact Riemannian manifold with 
$(g_{i j})_{i ,j = 1, \ldots, n}$ in local coordinates $ (f^1, \ldots, f^n) $.
Hermitian harmonic maps $f: M \longrightarrow N $ are defined as the solutions of 
the semi linear elliptic system 
  \begin{equation}
      \label{16}
      \gamma^{\alpha \bar{\beta}} ( \frac{\partial^2 f^i}{\partial z^\alpha \partial z^{\bar{\beta}}}                  			 +\Gamma^{i}_{j k}\frac{\partial f^j}{\partial z^\alpha}
      \frac{\partial f^k}{\partial z^{\bar{\beta}}} )=0, ~~i = 1, \ldots, n
 \end{equation}
This is the system  first studied by Jost and Yau \cite{JY}. As
discussed above, when
$M$ is not K\"ahler, the system (\ref{16}) is not in divergence form.  The method of Jost and Yau consists in 
studying the associated parabolic equation, 
  \begin{eqnarray}
    \label{17}
                                    f   &:& M \times [0, \infty) \longrightarrow N \\
        \frac{\partial f^i}{\partial t} &=& \gamma^{\alpha \bar \beta}
                                         (\frac{\partial^2 f^i}{\partial z^\alpha 
                                          \partial z^{\bar \beta}} +\Gamma^{i}_{j k}\frac{\partial f^j}
                                         {\partial z^\alpha}
                                          \frac{\partial f^k}{\partial x^{\bar \beta}} ) \\
                              f(z ,0)   &=& g(z)                       
    \end{eqnarray}
where $g: M \longrightarrow N$ is a continuous map. They show that a solution exists for $0 \leq t< \infty $,
under the assumption that $N$ has non positive sectional curvature and converges to a solution of (\ref{16}) 
under the geometric assumption of the following theorem:

\begin{theorem}[Jost-Yau]\label{JY}
Let M be a compact Hermitian manifold. Let $N$ be a compact Riemannian manifold of negative sectional
curvature. Let $g: M \longrightarrow N$ be continuous, and suppose
that $g$ is not homotopic to a map $g_0$ for which there is a nontrivial
parallel section of $g_0^{-1}TN$; for instance, assume that $g$ is not
homotopic to a map 
onto a closed geodesic of $N$. Then there is a  Hermitian harmonic map $f: M \longrightarrow N$ homotopic to 
$g$.
\end{theorem}     

In fact, an example in \cite{JY} shows that without this global
geometric assumption, a solution of the parabolic system need not
converge as $t\to \infty$, but may rather circle around $N$
forever. This is in contrast to the case of ordinary harmonic maps
where the variational structure forces a decay of the energy integral
along a solution of the parabolic flow which in turn implies that the
solution has to settle down asymptotically to a solution of the
elliptic system.   \\
As remarked above, when $M$ is K\"ahler, then something special
happens: The Hermitian harmonic map $f$ is simply an ordinary harmonic
map.

%%%%%%%%%%%%%%%%%%%%%%%%%%%%%%%%%%%%%%%%%%%%%%%%%%%%%%%%%%%%%%%%%%%%%%%%%%%%%%%%%%%%%%%%%%%%%%%%%%%%%%%%%%%%%%

\subsection{Affine harmonic maps}

As described above, on an affine manifold $M$ with metric tensor
$\gamma_{\alpha \beta}$, we can define an affinely 
invariant differential operator, 
$L := \gamma^{\alpha \beta} \frac{\partial^2}{\partial x^\alpha \partial x^\beta}$.  
A function 
$f: M \longrightarrow \mathbb R$ that satisfies $Lf = 0$ is called affine harmonic.
More generally, a map $f: M \longrightarrow N  $ where $N$ is a Riemannian manifold
with metric  $g_{i j}$ and Christoffel symbols $\Gamma^i_{j k}$ is called affine harmonic
if it satisfies 
\begin{equation}
      \label{18}
      \gamma^{\alpha \bar{\beta}} ( \frac{\partial^2 f^i}{\partial x^\alpha \partial x^{\bar{\beta}}}                  			 +\Gamma^{i}_{j k}\frac{\partial f^j}{\partial x^\alpha}
      \frac{\partial f^k}{\partial x^{\bar{\beta}}} )=0, ~~i = 1, \ldots, n
 \end{equation}
in local coordinates on $N$. In invariant notation (\ref{18}) can be written as 
\begin{equation} 
 \label{19}
   \gamma^{\alpha \beta} D_\alpha D_\beta f  = 0 
\end{equation}
where $D$ is the connection on the bundle $T^* M \otimes f^{-1}TN$ induced by the flat 
connection on $M$ and the Levi-Civit\'a connection on $N$. Jost and \c Sim\c sir 
obtained the following general existence result for affine harmonic maps, \cite{JS}.

\begin{theorem}[Jost - \c Sim\c sir]\label{JS}
Let $M$ be a compact  affine manifold, $N$ a compact Riemannian manifold of nonpositive sectional curvature. Let $g:M\to N$ be continuous, and suppose $g$ is not homotopic to a map $g_0:M\to N$ for which there is a nontrivial parallel section of $g_0^{-1}TN$. Then $g$ is homotopic to an affine harmonic map $f:M\to N$.
\end{theorem}

In fact, this result is stronger than the one stated in \cite{JS}; the
latter was formulated only for the special case of {\it K\"ahler}
affine manifolds in the sense of \cite{CY}. However, in the next
section, we shall describe the analytic scheme for showing existence
in such a way that it applies to any compact affine manifold $M$. 

Again, one may construct examples to show that the global topological
condition is needed in general, see \cite{JS}. Using the argument of Al'ber \cite{A}, one can also show that the affine harmonic map is unique in its homotopy class
under the assumptions of the above theorem. In fact, here, we also
need the global condition. 

%%%%%%%%%%%%%%%%%%%%%%%%%%%%%%%%%%%%%%%%%%%%%%%%%%%%%%%%%%%%%%%%%%%%%%%%%%%%%%%%%%%%%%%%%%%%%%%%%%%%%%%%%%%%%%%%%

\section{Analytic aspects of the existence scheme}
 
Consider the system (\ref{18}), (\ref{19})
\begin{equation}
 \label{20} 
    \frac{\partial f}{\partial t} = \gamma^{\alpha \beta}D_\alpha D_\beta f
\end{equation}

Linearizing  and using standard results about linear parabolic system,
which follow from the implicit function
theorem , it follows that (\ref{20}) has a solution for a short time interval $[0, \tau)$, and the 
interval of existence is open. Dealing with the global situation needs
the following steps which are harder. 
\begin{enumerate}
\item Showing the closedness of the existence interval, for which one needs the nonpositive sectional
curvature of the target manifold. 
\item Showing that the solution of (\ref{20}) converges to a non-divergence harmonic map as $t$ approaches
$\infty$, i.e., show that as $t$ approaches $\infty$, $\frac{\partial
  f}{\partial t}$ converges to 
$0$.
\end{enumerate}

In order to handle the first step one should show the local boundedness of the
energy density function 
$\eta(f) = \frac{1}{2} {\left\langle df, df \right\rangle}_{T^* M
  \otimes f^{-1} TN} $ where $df$ stands for the first derivatives of
$f$ w.r.t. the spatial variables $x$. 
For a detailed treatment of the procedure one may see \cite{JY} and \cite{JS}. Closedness of the
existence interval and thus the global existence follows from the regularity theory for 
 parabolic equations. In the following, we shall discuss the affine
 case; the complex case is analogous. Thus,  $x$ will now stand for
 affine coordinates. 

One first shows 
\begin{equation} 
 \label{21}
  \sup_{x\in M} g_{ij} \frac{\partial f^i}{\partial t} \frac{\partial
  f^j}{\partial t}
 \end{equation}
is nonincreasing in $t$. 
Next, $\eta(f)$ satisfies a linear
differential inequality, and we therefore obtain
\begin{equation}
   \label{22}
      \eta(f(x,t)) \leq c \sup_{t_0 \le \tau \leq t} \int_M \eta(f(.,\tau),
\end{equation}
for any $t_0 > 0$, see e.g. \cite{J1}, Section 3.3. Here and in the
sequel, $c$ stands for some constant that can be controlled by the
geometry of the manifolds involved, but which we do not make explicit here.\\
Next, using Jacobi field estimates \cite{J1} and the procedure in \cite{JS} one controls
the norm of $df$ with respect to the spatial variable $x$.
\begin{equation}
   \label{23}
      |df(x,t)|\leq c \left( \int_M \tilde{d}^2(f(.,\tau),f^0) \right)^{1/2} +c
\end{equation}   
where   $\tilde{d}(f(.,\tau),f^0) $is the homotopy distance between the initial map 
$f^0 = f(., 0)$ and the map $f(., t)$ at time $t$. It is defined as the length of the 
shortest geodesic from $f(x,t)$ to $f^0(x)$ in the homotopy class of curves 
determined by the homotopy between them.
Further computation leads to 
\begin{equation}
  \label{24}
|df(x,t)|\le c (1+t).
\end{equation}
Then, (\ref{21}) and (\ref{24}) yield $C^1$-bounds for the solution of (\ref{20}). 
In order to get $C^{2, \alpha} $ bounds, one may apply the regularity theory for solutions 
of linear parabolic equations by the standard bootstrapping argument.

For the second step of the proof one needs to show the convergence of the solution
$f(x,t)$ of (\ref{20}) to a non-divergence harmonic map at $\infty$. In this case, one needs
to require a topological non-triviality condition as expressed in the
Theorems \ref{JY}, \ref{JS} and also once more 
the nonpositive sectional curvature of the target manifold. 

We first choose a point $x_0 \in M$ where $\tilde{d}^2(f(y,t),f^0(y))$  attains its minimum
and apply the maximum principle on both the ball $B(x_0, R)$ of radius  $R$ and on its complement, to get
\begin{equation}
 \label{25}
\int_M \eta(f(.,t))\le c \sup_{y\in
  M}\tilde{d}(f(y,\tau),f^0(y)) +c. 
\end{equation}
Then (\ref{22}) gives the pointwise estimate
\begin{equation}
 \label{26}
   |df(x,t)|\le c (\sup_{y\in M}\tilde{d}(f(y,\tau),f^0(y)))^{1/2} + c.
\end{equation}
Therefore, for any $x_1,x_2 \in M$, denoting  the lift to universal covers  by  $\tilde{f}$ 
\begin{equation}
 \label{27}
    d(\tilde{f}(x_1,t),\tilde{f}(x_2,t))\leq c (\sup_{y \in M}\tilde{d}(f(y,\tau),f^0(y)))^{1/2} + c.
\end{equation}

The essential point of the proof then is to exclude that for
some sequence $t_n \to \infty$ for all $y\in M$,
\begin{equation}
  \label{28}
   \tilde{d}(f(y,t_n),f^0(y)) \longrightarrow \infty. 
\end{equation}
For the details, we  refer to \cite{JS}. For a family of solutions  $f(x, t, s) := f(x, t+s)$
depending on a parameter $s$, using (\ref{20}) 
\begin{equation}
  \label{29}
  \left( \gamma^{\delta \epsilon} \frac{\partial^2}{\partial x^\delta
    \partial x^\epsilon}- \frac{\partial}{\partial t} \right) \left(
  g_{ij} \frac{\partial f^i}{\partial s} \frac{\partial f^j}{\partial
    s} \right)\\
    =2\gamma^{\delta \epsilon}\left( g_{ij} \frac{\partial^2
    f^i}{\partial x^\delta \partial s} \frac{\partial^2 f^j}{\partial
    x^\epsilon \partial
    s} - \frac{1}{2} R_{ijkl}\frac{\partial f^i}{\partial
    s}\frac{\partial f^j}{\partial x^\delta}\frac{\partial
    f^k}{\partial s}\frac{\partial f^l}{\partial x^\epsilon}\right)  
 \end{equation}
one can conclude that, as $t$ tends to $\infty$, $\frac{\partial f(x, t)}{\partial t} $ 
converges to a parallel section $v(x)$ along $f_{\infty}$ which,
however,  is excluded in the assumptions of
Theorem \ref{JS}. Hence,  
\begin{equation}
  \label{30}
  \frac{\partial f(x,t)}{\partial t} \to 0 \text{ for } t \to \infty.
\end{equation}
This, together with the smooth convergence of $f(.,t_n)$ to
$f_\infty$,  shows that $f_\infty$ solves the elliptic system, i.e., it is affine
harmonic. \\
In fact, (\ref{29}) is also the key for the uniqueness of an affine
harmonic map in its homotopy class.  

%%%%%%%%%%%%%%%%%%%%%%%%%%%%%%%%%%%%%%%%%%%%%%%%%%%%%%%%%%%%%%%%%%%%%%%%%%%%%%%%%%%%%%%%%%%%%%%%%%%%%%%%%%%%%%%%%%%%

\section{Some possible future developments}
\begin{enumerate}
\item The theory of non-divergence harmonic maps can be investigated in a more general setting.
\item Dirichlet and Neumann boundary value problems for affine
  harmonic maps can be studied. 
In this case, the eternal circling of the solution is prevented by
the Dirichlet 
boundary values. Hence, here we do not need a global topological
condition. Of course, one now needs to prove  boundary regularity, but this problem can be solved by the methods of Jost and Yau \cite{JY},
or that of von Wahl \cite{W2}.
\item The method of Grunau and K\"uhnel \cite{GK} should be extended to show the existence of affine harmonic 
maps from a complete affine to a complete Riemannian manifold.
\item Most importantly, the results of Theorem \ref{JS} %and  of
%Corollary \ref{coro} 
should be applied to obtain rigidity results in affine 
differential geometry.
\end{enumerate}

%%%%%%%%%%%%%%%%%%%%%%%%%%%%%%%%%%%%%%%%%%%%%%%%%%%%%%%%%%%%%%%%%%%%%%%%%%%%%%%%%%%%%%%%%%%%%%%%%%%%%%%%%%%%%%%%%%%

\bibliographystyle{amsalpha}

\end{document}